\documentclass[11pt]{article}
\usepackage{graphicx}

\newtheorem{Theorem}{Theorem}

\def\qed{\hfill{\rm QED}}

\begin{document}

\pagestyle{empty}
\title{A Combinatorial Enumeration Approach For Measuring Anonymity}

\author{Jean-Charles Gr\'egoire\\
INRS--EMT\\
Montr\'eal, QC, Canada\\
gregoire@emt.inrs.ca\\
\and
 Ang\`{e}le M. Hamel\\
Physics and Computer Science\\
Wilfrid Laurier University\\
Waterloo, ON, Canada\\
ahamel@wlu.ca
}

\setcounter{secnumdepth}{2}

  \maketitle
\thispagestyle{empty}
    \begin{abstract}
    
A recent paper (Edman {\em et al.} \cite{ESY07}) has taken a combinatorial approach to measuring the anonymity of a threshold mix anonymous communications system.  Their paper looks at ways of matching individual messages sent to individual messages received, irrespective of user, and determines a measure of the anonymity provided by the system.  Here we extend this approach to include in the calculation information about how many messages were sent or received by a user and we define a new metric that can be computed exactly and efficiently using classical and elegant techniques from combinatorial enumeration.

    \end{abstract}

\section{Introduction}

Anonymity networks have evolved to address the problem of anonymous communication among users.  As internet technology becomes more prevelant in everyday life, questions of privacy and monitoring become more  important.  The anonymous network provides a means of communicating confidentially.  However, it is still vulnerable to attack.  One of the means of this is by attempting to match messages sent with messages received.  An exhaustive brute force attack is inefficient; statistical attacks are reasonably fast and reasonably effective.  Further, there is a need for a metric to measure the amount of anonymity that can be expected from a system.

A mix network, invented by Chaum \cite{C}, is a mechanism for anonymizing the correlation between  senders and receivers of messages.  Messages are sent into the mix where they are gathered, permuted and then delivered.  There are several mechanisms for doing this, including a threshold mix which takes in messages and holds them in a buffer until a predetermined threshold number of messages is reached and then it sends them.  The only possible attacks are based on observation of the input/output behaviour of the mix. We assume that it is possible for an adversary to see how many messages a user has sent and how many a user has received.

The main challenge for breaking anonymity in mix networks is determining whether some user Alice is communicating with some user Bob.  A secondary consideration can be the trajectory of a particular message, either who sent it or who received it.  Again focusing on Alice and Bob, we can either determine a {\em metric} which indicates the likelihood of correlating messages sent by Alice with messages received by Bob, or we can actually generate an {\em attack} which will attempt to break the system and reveal whom Alice is talking to (or who is talking to Bob).

Metrics  allow the user to make an informed choice of anonymity network.  They also allow an evaluation of how good an anonymity model is.  Historically, metrics have often considered either the perspective of an individual user or of an individual message.  A recent paper of Edman {\em et al.} expands the view to a system wide approach, but focuses on the traffic of individual messages.  Here we extend the approach to consider the traffic of sets of messages (in particular, messages sent by the same user, Alice, or received by the same user, Bob).

Section \ref{ourmetric} introduces our metric, reviews existing metrics and explores the differences ours manifests. Sections \ref{details} and \ref{extract} show how to calculate the metric. 
Section \ref{data} is a presentation and analysis of the data.
Section \ref{final} explores future work and delivers a conclusion.

\section{Our Metric}
\label{ourmetric}

 To establish our notation, suppose there are $k$ senders and $\ell$ receivers, the $i^{th}$ sender sends $s_i$ messages in a round and the $j^{th}$ receiver receives $r_j$ messages in a round, and that the total number of messages sent in a round is $n$, i.e. $s_1 + s_2 +\ldots +s_k =r_1+r_2 + \ldots +r_\ell = n$.  Then we want to know how many ways the $n$ messages could be divided up in this way.  These problems are often modeled, particularly in statistics literature, in terms of balls and urns.  In this language we want to know how many ways there are to deposit $n$ balls in $\ell$ urns, where there are $k$ different colours of balls: $s_1$ of one colour, $s_2$ of a second colour, etc., and each urn is to hold a particular number of balls: the first urn holds $r_1$ balls, the second urn holds $r_2$ balls, etc.

For example, suppose there are three messages, each labeled by $\alpha$, sent by $A_1$, three messages, each labeled by $\beta$, sent by $A_2$, and two messages, each labeled by $\gamma$, sent by $A_3$.  Suppose $B_1$ receives five messages and $B_2$ receives three messages.  Then by direct exhaustive, brute force enumeration of all the possibilities, there are nine different ways this could happen, where the first bracketing is the messages received by $B_1$ and the second bracketing is the messages received by $B_2$: $(\alpha, \alpha, \alpha, \beta, \beta) (\beta, \gamma,\gamma)$; $(\alpha,\alpha,\alpha,\beta,\gamma) (\beta, \beta, \gamma)$; $(\alpha,\alpha,\alpha,\gamma,\gamma) (\beta,\beta,\beta)$; $(\alpha,\alpha,\beta,\beta,\beta)(\alpha,\gamma,\gamma)$; $(\alpha, \alpha, \beta,\beta,\gamma) (\alpha, \beta,\gamma)$; $(\alpha,\beta,\beta,\beta,\gamma) (\alpha,\alpha,\gamma)$; $(\alpha, \beta, \beta, \gamma,\gamma) (\alpha, \alpha, \beta)$; $(\beta,\beta,\beta,\gamma,\gamma) (\alpha, \alpha,\alpha)$ $(\alpha, \alpha, \beta \gamma\gamma) (\alpha \beta\beta)$. 

In this type of system, the attacker can gain information by careful observation of the volume of messages originating or terminating at a user.  Consider at one extreme when $n$ messages are sent by Alice and $n$ messages are received by Bob and no messages are sent or received by any other users.  Then we know with certainty that the messages sent by Alice all went to Bob.  At the other extreme we have $n$ different senders each sending one message and $n$ different receivers each receiving one message.  In this case the number of possibilities for sender-receiver pairs is $n!$.   But even at intermediate stages, when some messages are sent by Alice and some messages are sent by others, and some messages are received by Bob and some messages are received by others, we can count the number of ways this could happen.
Counting this partitioning is actually a very old problem \cite{MacMahon} and can be solved in terms of a variety of generating functions called symmetric functions, as we will discuss in Sections \ref{details} and \ref{extract}. 

Our metric expresses the {\em degree of anonymity} as a ratio with the denominator representing the system with the most anonymity.  In this metric the most anonymity is provided by a system in which $n$ messages are sent but each sender sends exactly one message and each receiver receives exactly one message, as discussed above.  In this case there are $n!$ possibilities to match up a sender message with a receiver message.  We can informally define our metric as

\begin{equation}
\frac{\log(COUNT)}{\log(n!)},
\end{equation}

where $COUNT$ is the number of ways for $k$ senders and $\ell$ receivers to exchange $n$ messages if sender $i$ sends $s_i$ messages and receiver $j$ receives $r_j$ messages. 
We use the $\log$ here to have a compression of the scale for better
representation, and to avoid having numbers that are too large.
 In Section \ref{details} we will describe in detail how to calculate $COUNT$.

This metric is simple and straightforward to understand and calculate.  It is a system wide metric that measures the anonymity afforded by the system as a whole, rather than the anonymity afforded to a single user. 

We now review existing metrics.  Note that various perspectives are possible, e.g., the anonymity of an individual user, the anonymity of an individual message, or the anonymity of the system as a whole.  The anonymity metric seeks to distill into a single number the strength of the network with respect to protecting its users' anonymity.  This number is referred to as the degree of anonymity and was first proposed by Reiter and Rubin \cite{RR98}.  Their degree of anonymity requires a probability $p$ assigned to each potential sender and is defined as $1-p$ for each user.  A more systematic approach due to Berthold {\em et al.} \cite{BPS00} gives the degree of anonymity for a system of $N$ users as $A= \log_2 (N)$.  These metrics require estimates of properties of the system and can be imprecise.

The next step in metric evolution is the information theoretic (or entropy based) metrics. 
Serjantov and Danezis \cite{SD02} define a metric  $S=-\sum_{u=1}^{n} p_u \log_2 (p_u)$ where $n$ is the number of users, and $p_u$ is the probability that a user $u$ was the sender or the receiver in a message exchange. This metric is called the {\em effective anonymity set size} and it measures the entropy of the system. Recall that information entropy, as defined by Shannon, reflects the average information gained over a sequence of symbols, each having some probability. In this case, we have the probability that a specific user $u$ has sent a message. In the best case, this value will be equal to $\log_2(n)$, and grows with $n$.
In the worst case --- a user never sends, or is the only one to send --- it will be 0.

An improvement on this approach is the normalized metric of Diaz {\em et al.} \cite{DSCP02}, called the degree of anonymity, which is defined as
\[
\mbox{deg}=\frac{S}{S_{\mbox{max}}}= \frac{-\sum_{u=1}^{n} p_u \log_2 (p_u)}{\log_2 (n)} 
\]
where the term $S_{\mbox{max}}$ is the maximum entropy of the system which is $\log_2 (n)$.  This division
normalizes the result of Serjantov and Danezis, restricting the range to $[0,1]$ independent of $n$. 

There are drawbacks to these formulas. In this last case, it can be argued that, because of the normalization, it becomes easier to compare results but, at the same time, the number of users appears to become irrelevant. It also does not consider the users who have sent messages vs. the set of possible users: these metrics reflect a snapshot of the use of the system. Such a snapshot is also necessary to evaluate the probabilities required to compute the formula, but it is difficult to determine the degree of confidence we can have in such estimates, i.e. their quality.

A different approach is the combinatorial metric of Edman {\em et al.} \cite{ESY07}.
   They first define a bipartite graph, $G= (V_1, V_2, E)$, where $V_1$ is the set of sent messages and $V_2$ is the set of received messages.  There is an edge between two messages $s_i$ and $t_j$ if the sent message $s_i$ could be the same as the received message $t_j$.  Then this graph has an adjacency matrix, $A=(a_{i,j})_{n\times n}$, where the rows are indexed by the sent messages $V_1$ and the columns are indexed by the received messages $V_2$.  A perfect matching in a graph is a subset of edges such that every vertex is adjacent to exactly one edge in the subset.  In a bipartite graph this amounts to a pairing off of each vertex in set $V_1$ with exactly one vertex in set $V_2$.

In a bipartite graph it is well-known how to count the number of perfect matchings in the graph: they are counted by a mathematical function, the permanent, defined
\[
per(A) = \sum_{\sigma} \prod_{i=1}^n a_{i,\sigma(i)},
\]
where the sum is over all permutations, $\sigma$, of $n$ and $A=(a_{i,j})_{n\times n}$ is the adjacency matrix of the bipartite graph.  The reason the permanent works is as follows: every permutation selects an entry from each row and each column, so every pair consisting of a vertex in $V_1$ and a vertex in $V_2$ is represented exactly once.  If any of the selected entries is zero (i.e. there is no edge between those two vertices) then the product is zero and there is no perfect matching associated with that permutation. Conversely, if all entries are one then this permutation describes a perfect matching.

Edman {\em et al.} define a combinatorial degree of anonymity as follows:
\[
\mbox{deg} = \left\{ \begin{array}{cl}
0 & n=1\\
\frac{\log(per(A))}{\log(n!)} & n> 1\\
\end{array} \right.
\]
As with the degree of anonymity of Diaz {\em et al.} the measure reflects a ratio of the actual measurement over the ideal case.
The denominator is a reflection of the fact that the system providing the most anonymity is the one in which each sent message is potentially connected to each received message, i.e. the complete bipartite graph.  Then the $n\times n$ adjacency matrix is the all 1's matrix and the number of perfect matchings is equal to $n!$, the number of permutations of $n$.

Edman {\em et al.} then generalize their definition to matrices with  entries which are probabilities (doubly stochastic matrices).  In this model the probability in position $i, j$ is the probability that the edge between $s_i$ and $t_j$ is in a perfect matching.   Here, as in the unweighted case, they take a product of entries.  There are a number of concerns with this approach, some of which our approach corrects.

While the permanent counts perfect matchings in the unweighted case,
it is not clear which statistic is counted
in the weighted case, since
it is merely the sum of products of terms in the adjacency matrix. In the case of a 0--1 matrix the  permanent terms are the products of zeros and ones.  But the reason there is one term for each perfect matching is that this procedure is essentially a logical AND.  That is, a single zero will make the product of the entire set zero.   So while the entries are technically multiplied, they could just as easily be ANDed to the same effect.  In generalizing to the non 0-1 case it is not clear why multiplication should be the operation of choice to combine elements, nor what it counts.

Furthermore, this approach requires the calculation of probabilities for each edge and this in itself can be problematic.  Are these probabilities estimated? (with all the inherent issues of inaccuracy).  Are they calculated, say using an approach such as statistical disclosure?  If so, what is the  complexity of this approach and what are the quality of results it provides?

 Moreover, the perfect matching approach of Edman {\em et al.} considers a sent message, $s_i$, being matched to a received message and counts $s_i$ being matched to received message $t_u$ as different from $s_i$ being matched to $t_v$, even if $t_u$ and $t_v$ are received by the same user.  Certainly they are different messages, but if the goal is to determine who is communicating with whom, the important part is to determine that one of the many messages Alice sent is one of the many messages Bob received.

The recent work of Gierlichs {\em et al.} \cite{GTDPV08} refines the metric of Edman {\em et al.} to account for many messages sent and received by each user. 
To account for this Gierlichs {\em et al.} look at the equivalence class of perfect matchings.  This is actually the same situation as what we have already discussed. For example, if user $A_1$ sends $2$ messages and user $A_2$ sends $3$ messages, while user $B_1$ receives $2$ messages, user $B_2$ receives $2$ messages, and user $B_3$ receives one message, and each $A$ user could potentially communicate with each $B$ user, then there are $5!$ perfect matchings possible to pair up the messages sent with the messages received.

The authors denote the perfect matching by $M_C$ and each equivalence class by $[M_p]$ and cardinality $|[M_p]|=C_p$.  The total number of equivalence classes is some value $\Theta$ (determined by the problem).  In the example, then, there are $5$ equivalence classes with cardinalities $C_1=12, C_2=48, C_3=24, C_4=12, C_5=24$. 

\begin{eqnarray*}
\mbox{[$M_1$]}&=& [(A_1, B_1), (A_1, B_1), (A_2,B_2), (A_2, B_2),\\
& & (A_2, B_3)],\\
\mbox{[$M_2$]}&=& [(A_1,B_1), (A_1,B_2), (A_2,B_1), (A_2,B_1),\\
& & (A_2,B_2), (A_2,B_3)],\\
\mbox{[$M_3$]}&=& [(A_1,B_1),(A_1,B_3), (A_2,B_1),(A_2,B_2),\\
& &(A_2,B_2)],\\
\mbox{[$M_4$]}&=& [(A_1,B_2),(A_1,B_2),(A_2,B_1),(A_2,B_1),\\
& &(A_2,B_3)],\\ 
\mbox{[$M_5$]}&=& [(A_1,B_2),(A_1,B_3),(A_2,B_1),(A_2,B_1),\\
& &(A_2,B_1), (A_2,B_2)]\\
\end{eqnarray*}

where we have used the notation $(A_i,B_j)$ to mean a pair consisting of an element from $A_i$ and an element from $B_j$, so, for example, $(A_2,B_1)$ is the set of all pairs $ (a_{2}^{1}, b_{1}^{1}), (a_{2}^{1}, b_{1}^{2}), (a_{2}^{2}, b_{1}^{1}),$ $ (a_{2}^{2}, b_{1}^{2}), (a_{2}^{3}, b_{1}^{1}), (a_{2}^{3}, b_{1}^{2})$ where $a_{2}^{i}$ means the $i$th message sent by user $A_2$.

The authors define the system's anonymity level, $d^*(A)$, as 
\[
\frac{-\sum_{p=1}^{\Theta} Pr(M_C\in [M_{p}])\cdot \log(Pr(M_C\in[M_{p}]))}{\log(n!)} 
\]
if $n>1$, and as $0$ if $n=1$, where $Pr(M_C\in [M_p])= \frac{C_p}{per(A)}$.

 It is worth noting that Franz {\em et al.} \cite{Franz} also take a counting approach, although  they do not do the full generality of sender/receiver patterns.  In one instance they look at all possible combinations of senders for a given set of messages.  In another instance they count senders sending  various combinations of messages but do not consider receivers receiving several messages.   While some ideas are similar to ours and could be expanded further using enumerative techniques, they do not have the full generality of our approach.  In fact, classical enumerative methods used by Franz {\em et al.} and Gierlichs {\em et al.} indeed work to count when either the senders or the receivers are fixed at sending or receiving one message each; to do the full generality of both receivers and senders together dealing in multiple messages one needs symmetric functions, as discussed below.

 Like Gierlichs {\em et al.} we suppose that senders and receivers send many messages, we have a combinatorial metric like Edman {\em et al.} and ask ``how many ways for these senders to send to these receivers?". This differs from Gierlichs {\em et al.} who use an entropy based metric, asking instead, ``what's the probability that this perfect matching is right one?"  Their approach requires the calculation of two parameters: equivalence classes and cardinality (see Appendix \ref{appendixB} for a discussion of a way of calculating cardinality).  To do this they provide a divide-and-conquer algorithm that they note becomes rather expensive for large $n$.  Indeed in their conclusions they suggest that a more efficient algorithm remains an open problem.   Our method essentially calculates the size of the equivalence classes; however, it does so without explicitly enumerating them, thus the approach is extremely fast and streamlined.  It provides a rapid but accurate measure of the anonymity of the system.  As we discuss in Section \ref{data}, numerous trends can be discerned and this provides an interesting focus for future work.

\section{Calculating our Metric}
\label{details}

We now turn our attention to determining $COUNT$ as defined in Section \ref{ourmetric}.  The calculation of it is straightforward; however, it requires some ``heavy machinery'' from combinatorial enumeration, namely generating functions and symmetric functions (which are a special type of generating function).
 We briefly review generating functions before discussing the appropriate one for this particular problem.  Excellent introductions to generating functions can be found in \cite{GouldenJackson}, \cite{Stanley} or \cite{Wilf}. 

A generating function is a sum of powers of $x$ where the coefficient of $x^i$ counts how many items of size $i$ there are.  
In a sense the powers of $x$ are merely placeholders, with the $i^{th}$ power holding the place for items of size $i$, and the $x$'s are not expected to be evaluated.  For example, if there are four ways of having two messages delivered, three ways of having one message delivered, and one way of having no messages delivered, then the generating function is $1+ 3x + 4x^2$. 
For a counting problem a generating function is set up that models the problem and then the required coefficient is extracted.  
  The notation $[x^i]$ means ``the coefficient of $x^i$.''
Thus in the example, $[x^2](1+3x+4x^2)$ will give us the value $4$.

Generating functions have the advantage that they encode all the enumerative information and they can easily be manipulated, e.g. multiplied together.  The extraction of a coefficient can prove to be a challenge sometimes if a direct formula for it is not easily obtainable.  In this case a symbolic computation program such as Maple can be an important tool.

Turning to our specific problem, the generating function allows us, given the number of messages sent and received by various users, to determine exactly the number of ways this could take place.  If there are a lot of ways for this to take place then the system does not leak much information and remains relatively anonymous.  If there are only a few ways for this to take place then the system is leaking a lot of information.

We define our degree of anonymity precisely as follows.  As mentioned above, we take a ratio with the denominator representing the system with the most anonymity, i.e a system in which each user sends or receives a single message.   In this system our approach is no better than counting perfect matchings and there are $n!$ possibilities.  The numerator is the number of ways the $n$ messages could be divided up.  This is the coefficient of $x_{1}^{s_{1}}x_{2}^{s_{2}}\ldots x_{k}^{s_{k}}$ in the generating function, $GF$, for the number of ways of $n$ messages being received such that the $i$th user receives $r_i$ messages, for $1\leq i\leq \ell$.
Thus the degree is
\begin{equation}
deg_{A} = \frac{\log([x_{1}^{s_{1}}x_{2}^{s_{2}}\ldots x_{k}^{s_{k}}] GF)}{\log(n!)}.
\end{equation}

Now we consider the form of the generating function, $GF$, for this problem. This generating function is a special type of function called a {\em symmetric function}. First, a number of further definitions are required.  A  symmetric function $f(x)$ in variables $x_1, x_2, \ldots x_k$ is a function such that a permutation of the variables does not change the value of the function, i.e. $f(x_{\sigma(1)}, x_{\sigma(2)},\ldots, x_{\sigma(k)})=f(x_1,x_2,\ldots, x_k)$.  Then the {\em homogeneous symmetric function} of degree $m$, $h_m({\bf x})$, is the sum of all homogeneous terms in $x=x_1, x_2,\ldots, x_k$, i.e.
\begin{eqnarray*} \small
h_1(x_1, x_2, x_3)&=& x_1 + x_2 +x_3\\
h_2(x_1, x_2, x_3) &=& x_{1}^{2} + x_{2}^{2} + x_{3}^{2} + x_1x_2\\
& & +\; x_1 x_3 + x_2 x_3\\
h_3(x_1, x_2, x_3) &=& x_{1}^{3} + x_{2}^{3} + x_{3}^{3} + x_{1}^{2} x_{2} + x_{1}^{2}x_3 \\
&& +\; x_1 x_{2}^{2}   + x_{2}^{2} x_3 + x_{1} x_{3}^{2}+ x_2 x_{3}^{2}\\
& & +\; x_1 x_2 x_3\\
\end{eqnarray*} 

The homogeneous symmetric functions can also be defined for infinite sets of variables.  Furthermore, for  $\lambda=\lambda_1, \lambda_2,\ldots, \lambda_m$ a partition of $n$, where $\lambda_1+ \lambda_2+\ldots +\lambda_m = n$ (i.e.\ a nondecreasing sequence of nonnegative integers that sum to $n$), then $h_\lambda$ is defined as the product $h_{\lambda_{1}} h_{\lambda_{2}} \ldots h_{\lambda_{m}}$. 


\begin{Theorem}
Given that $s_i$ messages are sent by each sender, $1\leq i\leq k$, and that $r_i$ messages are received by each receiver, $1\leq i\leq\ell$, in a round, the number of ways this could happen is
\[
[x_{1}^{s_{1}} x_{2}^{s_{2}}\ldots x_{k}^{s_{k}}] h_{r_{1}} ({\bf x}) h_{r_{2}} ({\bf x}) \ldots h_{r_{\ell}} ({\bf x})
\]
where ${\bf x}$ is $x_1, x_2, \ldots, x_k$ and $h_q(x)$ is the homogeneous symmetric function of degree $q$.
\label{thm1}
\end{Theorem}

Proof:
The term $h_m(x_1, x_2\ldots, x_k)$ counts the number of different ways $m$ elements could be received where the elements are drawn from $1, 2,\ldots, k$ (e.g. if there were three elements received and two possible kinds of elements, then this is $h_3(x_1, x_2)= x_{1}^{3} + x_{2}^{3} + x_{1}^{2} x_2 + x_1 x_{2}^{2}$).

The product $h_{r_{1}}({\bf x}) h_{r_{2}}({\bf x})\ldots h_{r_{\ell}}({\bf x})$ is, by the product lemma in enumerative combinatorics \cite{GouldenJackson}[pp36-37], the generating function for the number of ways of one user receiving $r_1$ elements, a second user receiving $r_2$ elements, etc, simultaneously.

The term $[x_{i}^{s_{i}}]$ denotes the coefficient of $x_{i}^{s_{i}}$ in the expression.  This counts the number of ways $s_i$ $i$'s could be sent.

The entire expression in the statement of the theorem thus counts in general the number of ways $s_i$ $i$'s could be sent and $r_i$ $i$'s could be received.

\qed

Thus the generating function we require is a symmetric function and  we can define our degree of anonymity to be
\begin{equation}
deg_A=  \frac{\log([x_{1}^{s_{1}}x_{2}^{s_{2}}\ldots x_{k}^{s_{k}}] h_{r_{1}} ({\bf x}) h_{r_{2}} ({\bf x}) \ldots h_{r_{\ell}} ({\bf x})
)}{\log(n!)}.
\label{degA}
\end{equation}
where there are $k$ senders and $\ell$ receivers, the $i^{th}$ sender sends $s_i$ messages in a round and the $j^{th}$ receiver receives $r_j$ messages in a round, and that the total number of messages sent in a round is $n$, i.e. $s_1 + s_2 +\ldots +s_k =r_1+r_2 + \ldots + r_\ell = n$.

Consider an example. Recall that earlier  we showed the example of an input/output round of three messages sent by $A_1$, three messages sent by $A_2$, two messages sent by $A_3$, five messages received by $B_1$, and three messages received by $B_2$.  In our generating function terms this means that we need two complete generating functions: $h_5$ for $B_1$ and $h_3$ for $B_2$ (since $B_1$ receives five messages and $B_2$ receives three messages).  Since there are three users sending messages, the number of variables for each generating function is limited to three.   Since we know that $A_1$ sends three messages, $A_2$ sends three messages, and $A_3$ sends two messages, we require the coefficient of $x^{3}_{1}x_{2}^{3}x_{3}^{2}$.  Specifically,

\begin{eqnarray*}
h_5(x_1, x_2,x_3) & = & x_{1}^{5}+x_{2}^{5}+x_{3}^{5}+x_{1}^{4}x_{2}+x_{1}^{4}x_{3}\\
& &+x_{2}^{4}x_{1}
+ x_{2}^{4}x_{3} +x_{3}^{4}x_{1}+x_{3}^{4}x_{2}\\
& & +x_{1}^{3}x_{2}^{2}+x_{1}^{3}x_{3}^{2} +x_{2}^{3}x_{1}^{2}+x_{2}^{3}x_{3}^{2}\\
& &+x_{3}^{3}x_{1}^{2}+x_{3}^{3}x_{2}^{2}
+x_{1}^{3}x_{2}x_{3}\\
 & & +x_{2}^{3}x_{1}x_{3}+x_{3}^{3}x_{2}x_{1}+x_{1}^{2}x_{2}^{2}x_{3}\\
& &+x_{1}^{2}x_{3}^{2}x_{2}+x_{2}^{2}x_{3}^{2}x_{1}
\end{eqnarray*}

and

\begin{eqnarray*}
h_{3}(x_1,x_2,x_3) &= & x_{1}^{3}+x_{2}^{3}+x_{3}^{3}+x_{1}^{2}x_{2}+x_{1}^{2}x_{3}\\
& &+x_{2}^{2}x_{1} +x_{2}^{2}x_{3}+x_{3}^{2}x_{1}+x_{3}^{2}x_{2} \\
& &+x_1x_2x_3
\end{eqnarray*}

We can multiply these two generating functions together and collect terms (admittedly a slow process, but we will improve on it in Section \ref{extract}).  This approach shows that  $h_5(x_1, x_2, x_3) h_3 (x_1, x_2, x_3)$ has nine terms of the form $x_{1}^{3}x_{2}^{3}x_{3}^{2}$ formed from the following products: $(x_{1}^{3}x_{2}^{2}) (x_{2}x_{3}^{2})$;
 $(x_{1}^{3}x_{2}x_{3}) (x_{2}^{2}x_{3})$; 
$(x_{1}^{3}x_{3}^{2}) (x_{2}^{3})$;
 $(x_{2}^{3}x_{1}^{2})(x_{1} x_{3}^{2})$;
$(x_{1}^{2}x_{2}^{2}x_{3})(x_{1}x_{2}x_{3})$;
$(x_{2}^{3}x_{1}x_{3})(x_{1}^{2}x_3)$; 
 $(x_1x_{2}^{2}x_{3}^{2})(x_{1}^{2} x_2)$;
 $(x_{2}^{3}x_{3}^{2})(x_{1}^{3})$
 $(x_{1}^{2}x_{2}x_{3}^{2})(x_{1}x_{2}^{2})$; .\\
Thus $[x_{1}^{3} x_{2}^{3} x_{3}^{2}] h_5(x_1, x_2, x_3) h_3 (x_1, x_2, x_3)=9$. Note that these can be matched exactly with the $\alpha,\; \beta,\; \gamma$ terms obtained earlier in the section through direct enumeration of the various possibilities.

Also compare this approach with the permanent-based approach that considers all possible combinations of messages sent and received.  Since there are eight messages involved, there are $8!=40320$ ways to send them, a substantially larger number of possibilities.  Our degree of anonymity is $\log 9/\log 40320 =   0.954/4.605=0.207$ whereas the degree of anonymity of Edman {\em et al.} is $1$.

\section{Extracting the Coefficient}
\label{extract}

Recall from Section \ref{details} that the generating function for the problem is the homogeneous symmetric function $h_\lambda$ and that in order to evaluate our degree of anonymity, $deg_A$, in equation (\ref{degA}) we need to extract the coefficient.  This section explores the theoretical basis for this extraction and explains the calculation that needs to be made.  

Symmetric functions form a \emph{graded ring}. The most natural basis for this ring is the set of {\em monomial} symmetric functions.  The monomial symmetric functions, $m_\lambda$, are defined as $m_\lambda({\bf x})= \sum_{\alpha} {\bf x}^{\alpha}$ where the sum ranges over all distinct permutations $\alpha=(\alpha_1, \alpha_2, \ldots, \alpha_n)$ of the entries of the partition $\lambda=(\lambda_1, \lambda_2,\ldots, \lambda_n)$. For example, $m_{3,1,1}(x_1,x_2,x_3,x_4) =x_{1}^{3} x_{2} x_{3} + x_1^{3} x_2 x_4 + x_{1}^{3} x_3 x_4 + x_{2}^{3} x_1 x_3 + x_{2}^{3} x_1 x_4 + x_{2}^{3} x_3 x_4 + x_{3}^{3} x_1 x_2 + x_{3}^{3} x_1 x_4 + x_{3}^{3} x_2 x_4 + x_{4}^{3} x_1 x_2 + x_{4}^{3} x_1 x_3 + x_{4}^{3} x_2 x_3$.  As a basis, then, we can write any symmetric function, $f({\bf x})$, as $\sum_\lambda f_\lambda m_\lambda ({\bf x})$ where the sum is over all partitions of $n$.   

There is also a very natural scalar product defined on this ring. It is defined such that $<m_\lambda ({\bf x}), h_\mu({\bf x})>= \delta_{\lambda,\mu}$ where $\delta_{\lambda, \mu}$ equals $1$ if $\lambda=\mu$ and $0$ otherwise. Note in particular that, although the monomial symmetric functions, and indeed the homogeneous symmetric functions, are both bases for the ring of symmetric functions, neither is an orthonormal basis with this scalar product. This scalar product however allows us to extract the coefficients of a symmetric function $f({\bf x})$.  Suppose we want $[x_{1}^{\mu_{1}} x_{2}^{\mu_{2}}, \ldots x_{k}^{\mu_{k}}]$.  Then 
\begin{eqnarray}
<f({\bf x}), h_\mu({\bf x})> & =& < \sum_\lambda f_\lambda m_\lambda ({\bf x}), h_\mu ({\bf x})>\\ 
 &=& \sum_\lambda f_\lambda <m_\lambda ({\bf x}), h_\mu ({\bf x})>\\
& =& f_\mu. 
\label{scalcoeff}
\end{eqnarray}

Thus to compute the number of ways $k$ senders could send $s_1, s_2, \ldots, s_k$ messages and $\ell$ receivers could receive $r_1,r_2, \ldots, r_\ell$ messages, we calculate the scalar product $<h_{s_{1}}*h_{s_{2}}*\ldots *h_{s_{k}}, h_{r_{1}}*h_{r_{2}}*\ldots* h_{r_{\ell}}>$. 

The step-by-step justification for the procedure to calculate the number of ways $k$ senders and $\ell$ receivers send $n$ messages in a round such that the $i^{th}$ sender sends $s_i$ messages and the $j^{th}$ receiver receives $r_j$ messages, is as follows:

\begin{enumerate}
\item By Theorem \ref{thm1} this number can be represented by $[x_{1}^{s_{1}} x_{2}^{s_{2}}\ldots x_{k}^{s_{k}}] h_{r_{1}} ({\bf x}) h_{r_{2}} ({\bf x}) \ldots h_{r_{\ell}} ({\bf x})$.
\item By equation (\ref{scalcoeff}) the coefficient for $x_{1}^{s_{1}} x_{2}^{s_{2}}\ldots x_{k}^{s_{k}}$ is equal to the scalar product of the generating function representing senders and the generating function representing receivers.
\item By the proof of Theorem \ref{thm1} and the comments before Theorem \ref{thm1}, the complete symmetric function, $h_s$, is the generating function for a sender sending $s$ messages, and $h_{s_{1}}h_{s_{2}}\ldots h_{s_{k}}$ is the generating function for $k$ senders with the $i^{th}$ sender sending $s_i$ messages.  Similarly for receivers.
\end{enumerate}

Now that we can extract the coefficient via the scalar product, we can calculate our degree of anonymity, $deg_A$.  However, with the exception of a few special cases, obtaining a closed form expression for the scalar product is difficult.  The alternative is to use a symbolic computation package, such as Maple, to calculate the scalar product.  In the next section we outline the results we obtained using such a program. The computations presented here were carried out using Maple 8 and the symmetric functions package, SF, written by John Stembridge \cite{Stembridge}. All of the calculations mentioned ran in a few seconds or less.

\section{Data Analysis}
\label{data}

 We have already discussed the extremes of the metric (i.e. when Alice sends all messages, or when each user sends exactly one) and have discussed how it discerns between cases better than the metric of Edman {\em et al.}   In this section we conduct a number of experiments on the metric and discover a number of patterns and trends.  In particular we explore what the metric looks like, some interesting features of it, and answers to some interesting questions.

 The two figures, Figure 1 and Figure 2, illustrate the
 behaviour of the metric from more favourable (anonymous) to less favourable
 situations, based on two different scenarios where Alice's communications
become predominant.   The first considers the case where the number
of messages sent by Alice increases but the number of messages stays
the same.  The second considers the case where the number of messages
sent by Alice increases and the total number of messages also
increases.  These both show that
 the metric tends to linear as the values increase away from perfect
 anonymity.

\begin {figure}[ht]
\centering
  \includegraphics[angle=-90,scale = 0.3]{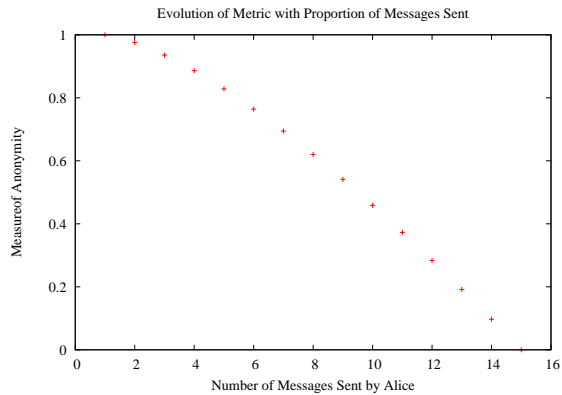}
\caption {For a total 15 messages sent to distinct receivers, evolution as Alice sends growing numbers of these messages.}
\label {fig:1}
\end{figure}

\begin {figure}[ht]
\centering
  \includegraphics[angle=-90,scale = 0.3]{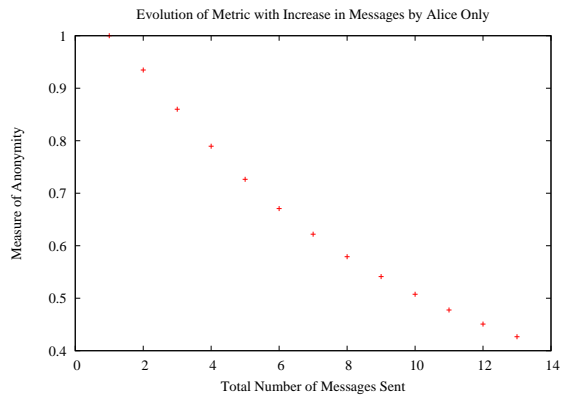}
\caption {For seven distinct senders, including Alice, and one receiver per message, evolution as Alice sends from one to 13 messages and all others send one.}
\label {fig:2}
\end{figure}

Now we turn to features. This metric has many interesting features, as illustrated by the following example for $n=7$. We have calculated the coefficients for all sender-receiver combinations for $n=7$ messages and the results are shown in a table in the Appendix. Note that in the table we have taken advantage of symmetry, i.e. $s_1,s_2\ldots, s_k$ senders and $r_1,r_2,\ldots, r_\ell$ receivers gives the same coefficient as $r_1, r_2,\ldots, r_k$ senders and $s_1,s_2,\ldots,s_k$ receivers.  
We have worked with $n=7$ because the number of cases to consider is tractable and the coefficients obtained are small enough to be meaningful; however, the calculations also run in seconds on larger values of $n$.  
The table shows clearly that, as one would expect, a lot of senders sending a few messages each (or a lot of receivers, receiving few messages each) results in the most anonymity.  A closer examination of the table reveals a number of interesting facts.

First, coalescing a sender of a single message into another sender  (e.g. going from $1,1,1,1,1,1,1$ to $1,1,1,1,1,2$) cuts the coefficient by a factor proportional to the new sender's  number of messages, a dramatic reduction in the anonymity.  More precisely, going from $1, k$ to $k+1$ divides the coefficient by $k+1$.  The same applies to receivers.  The coefficient is cut by a smaller amount because the log smooths the function.

 Second, there are counterintuitive instances as well where more senders (resp. receivers) and fewer messages is not superior.  For example, $1,1,1,1,1,1,1; 1,2,2,2$ has a coefficient of $630$ and degree $0.756$ and $1,1,1,1,1,2; 1,1,1,2,2$ has a coefficient of $690$ and coefficient $0.767$, yet the first sender-receiver pair has more senders.  Of course the second has more receivers. A further interesting case is $1,3,3; 2,2,3$ with coefficient $19$ and degree 0.345, and $1,3,3; 1, 1,1,4$ with coefficient $20$ and degree $0.351$, and thus a barely discernable difference between the coefficients and degree, yet the second group has more receivers.  It can be difficult to strike a hard and fast rule about which is better.

Third, it is actually relatively easy to rank the partitions from most anonymous to least.  

 \begin{center}
\begin{tabular}{|l|l|l|}
\hline
Messages & Coeff& Deg\\ \hline
1,1,1,1,1,1,1 & 5040 & 1 \\ \hline
 1,1,1,1,1,2 & 2520 & 0.919 \\ \hline
 1,1,1,2,2 & 1260 & 0.832\\ \hline
 1,1,1,1,3 & 840 & 0.790 \\ \hline
 1,2,2,2 & 630 & 0.756\\ \hline
 1,1,2,3 & 420 & 0.708 \\ \hline
 1,1,1,4 & 210& 0.627\\ \hline
 2,2,3 & 210& 0.627 \\ \hline
 1,3,3 & 140 & 0.580\\ \hline
 1,2,4 & 105& 0.546\\ \hline
 1,1,5 & 42 &0.438\\ \hline
 3,4 & 35& 0.417 \\ \hline
 2,5& 21& 0.357\\ \hline
 1,6 & 7& 0.228 \\ \hline
 7 & 1& 0 \\ \hline
\end{tabular}
\end{center}

In comparing the 15 possible sender partitions with a receiver set $1,1,1,1,1,1,1,$ we can order, from most anonymous, to least anonymous as shown in the table above.
This trend holds up with other receiver sets as well, showing that, in general, more senders is better, although, as noted in the second point, when the receiver set also varies, there is some variation from this rule.

\begin {figure}[ht]
\centering
  \includegraphics[angle=-90,scale = 0.3]{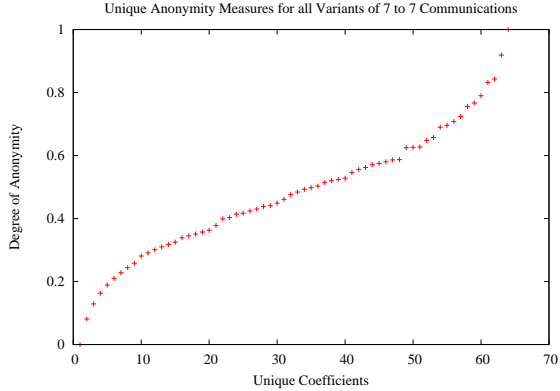}
\caption {Ordered Values of the Metric Generated from All Permutations of Seven}
\label {fig:6}
\end{figure}

Finally, we can also represent the $n=7$ information  graphically in Figure \ref{fig:6}.  This plots the different coefficients we may have and shows that, except near the extreme values, the changes in the degree tend to be linear in terms of the variations in the coefficients.

Another interesting feature for any $n$ occurs when the number of messages sent by one party is equal to or greater than the total messages sent by all other parties.  In this case the coefficient stops increasing, although, of course, the degree decreases.  For example (note we used the notation $1^6$ to mean $1,1,1,1,1,1$),

\begin{center}

\begin{tabular}{|l|r|l|}
\hline
Messages & Coeff & Deg\\ \hline
1, $1^6$; 1, $1^6$  &   5040 & 1\\ \hline
2, $1^6$; 2, $1^6$ &   10,440 & 0.873\\ \hline
3, $1^6$; 3, $1^6$  &  12,840& 0.739\\ \hline
4, $1^6$; 4, $1^6$ &    13,290& 0.629 \\ \hline
5, $1^6$; 5, $1^6$ &     13,326& 0.543\\ \hline
6, $1^6$; 6, $1^6$  &     13,327& 0.475\\ \hline
7, $1^6$; 7, $1^6$ &     13,327&0.421 \\ \hline
8, $1^6$; 8, $1^6$  &     13,327&0.377 \\ \hline
9, $1^6$; 9, $1^6$  &     13,327&0.340 \\ \hline
10, $1^6$; 10, $1^6$  &    13,327&0.310 \\ \hline
\end{tabular}
\end{center}

\begin {figure}[ht]
\centering
  \includegraphics[angle=-90,scale = 0.3]{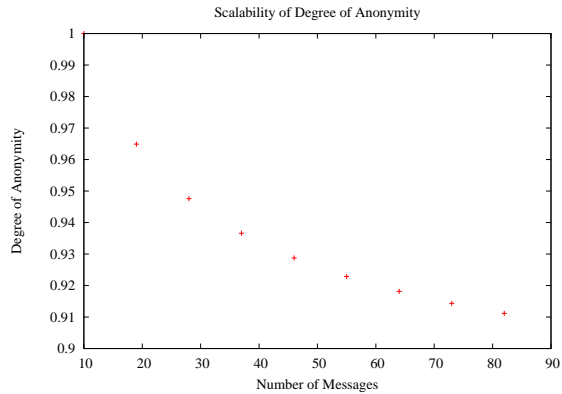}
\caption {Scalability of the metric}
\label {fig:5}
\end{figure}

In terms of answering an interesting question, we can derive some information on the importance of the size of the mix on the metric.
As any degree metric is a ratio, it does not explicitly take into account the number of users in the system. It is then a legitimate question to see how, for a given sender/receiver pattern, the value evolves as the size of the mix increases.
Figure \ref{fig:5} presents such an experiment, where the ratio of messages sent from Alice to Bob remains the same as the number of users increases. 
More precisely, the figure shows the effect on degree of anonymity as the
number of messages increases in the case when the ratio of Alice's
messages sent (and Bob's messages received) to the total number of
messages sent stays the same.  In each  instance $k$ messages are
sent: Alice sends $p$ messages, and the $k-p$ other users each send
one message, while Bob receives $p$ messages and the other $k-p$ users
each receive one message (for $p=1,2,3,4,5,6,7,8,9$ and $k=9p+1$).

\section{Conclusion and Future Work}
\label{final}
We have developed an elegant and easy way to calculate a metric for the degree of anonymity of an anonymous communication system.  Our metric uses techniques from classical enumeration to count without actually calculating the various possibilities of senders and receivers exchanging different combinations of messages.  It nicely and naturally extends the metric of Edman {\em et al.} \cite{ESY07} who introduced the combinatorial approach to this area.
Our metric is straightforward to calculate employing Maple and using it we are able to produce data that highlight a number of significant trends such as we produced in Section \ref{data}.

Future work will focus on practical uses of the metric.
One possible direction is to turn our metric into an attack.  At the moment it is a measure of the degree of anonymity a user could expect, but it may be possible to further exploit the knowledge of the partition sizes to break the anonymity network.

The ``mostly linear'' behaviour we have mentioned in Section \ref{data} deserves further study. While we feel linearity is an interesting property for a metric to have, Figure \ref{fig:6} shows that it does not appear in the extreme cases of full or zero anonymity. Further investigation is required to explore the limits of this approximation with large sizes, and how it can be best exploited practically.

A further useful practical outcome would be a recommendation to Alice on what she should do at each stage.  For example, would sending another message to Bob increase the likelihood of detection or not?  Is she advised to send a message to someone else?  Should she wait until the next round?  It may be possible to provide guidance to her using knowledge of the partition sizes.

Finally, future work could include more analysis of the data: this metric is a fast, simple tool for calculating the anonymity of the system;  as we showed in Section \ref{data} it allows a number of interesting features to be detected.  More analysis will yield more patterns.

\vfill

\appendix
\section{Example}

This is the table of all possible sender-receiver pairs for $n=7$ messages.  The table includes the symmetric function coefficient as well as the degree of anonymity. Note that the notation $1^7$ means $1,1,1,1,1,1,1$.

\begin{tabular}{|l|l|l|}
\hline
Messages & Coeff& Deg\\ \hline
$1^7$;$1^7$ & 5040 & 1 \\ \hline
$1^7$; $1^5$,2 & 2520 & 0.919 \\ \hline
$1^7$; $1^3$,2,2 & 1260 & 0.832\\ \hline
$1^7$; 1,1,1,1,3 & 840 & 0.790 \\ \hline
$1^7$; 1,2,2,2 & 630 & 0.756\\ \hline
$1^7$; 1,1,2,3 & 420 & 0.708 \\ \hline
$1^7$; 1,1,1,4 & 210& 0.627\\ \hline
$1^7$; 2,2,3 & 210& 0.627 \\ \hline
$1^7$; 1,3,3 & 140 & 0.580\\ \hline
$1^7$; 1,2,4 & 105& 0.546\\ \hline
$1^7$; 1,1,5 & 42 &0.438\\ \hline
$1^7$; 3,4 & 35& 0.417 \\ \hline
$1^7$; 2,5& 21& 0.357\\ \hline
$1^7$; 1,6 & 7& 0.228 \\ \hline
$1^7$; 7 & 1& 0 \\ \hline
$1^5$,2; $1^5,2$ & 1320& 0.843\\ \hline
$1^5$,2; 1,1,1,2,2 & 690 &0.767\\ \hline
$1^5$,2; 1,1,1,1,3 & 480 &0.724\\ \hline
$1^5$,2; 1,2,2,2 & 360 &0.690\\ \hline
$1^5$,2; 1,1,2,3 & 250 &0.648\\ \hline
$1^5$,2; 1,1,1,4 & 135 &0.575\\ \hline
$1^5$,2; 2,2,3 & 130 & 0.571 \\ \hline
$1^5$,2; 1,3,3 & 90 &0.528\\ \hline
$1^5$,2; 1,2,4 & 70 &0.498\\ \hline
$1^5$,2; 1,1,5 & 31 &0.403\\ \hline
$1^5$,2; 3,4 & 25 &0.378\\ \hline
$1^5$,2; 2,5 & 16 &0.325\\ \hline
$1^5$,2; 1,6 & 6 & 0.210\\ \hline
$1^5$,2; 7 & 1& 0\\ \hline
1,1,1,2,2; 1,1,1,2,2 & 378& 0.696 \\ \hline
1,1,1,2,2; 1,1,1,3 & 270 &0.657\\ \hline
1,1,1,2,2; 1,2,2,2& 207 & 0.626\\ \hline
1,1,1,2,2; 1,1,2,3 & 148& 0.586\\ \hline
1,1,1,2,2; 1,1,1,4 & 84 &0.520\\ \hline
1,1,1,2,2; 2,2,3 & 81 & 0.515 \\ \hline
1,1,1,2,2; 1,3,3 & 58 &0.476\\ \hline
1,1,1,2,2; 1,2,4 & 46 &0.449\\ \hline
1,1,1,2,2; 1,1,5 & 22& 0.363 \\ \hline
1,1,1,2,2; 3,4 & 18 &0.339\\ \hline
1,1,1,2,2; 2,5 & 12 &0.291\\ \hline

\end{tabular}

\begin{tabular}{|l|l|l|}
\hline
Messages & Coeff&Deg\\ \hline
1,1,1,2,2; 1,6 & 5 &0.189\\ \hline
1,1,1,2,2; 7& 1 & 0 \\ \hline
1,1,1,1,3; 1,1,1,1,3 & 208 &0.626\\ \hline
1,1,1,1,3; 1,2,2,2 & 150 & 0.588 \\ \hline
1,1,1,1,3; 1,1,2,3 & 114 &0.556\\ \hline
1,1,1,1,3; 1,1,1,4 & 73 & 0.503 \\ \hline
1,1,1,1,3; 2,2,3 & 62&0.484 \\ \hline
1,1,1,1,3; 1,3,3 & 46 &0.449\\ \hline
1,1,1,1,3; 1,2,4 & 39 &0.430\\ \hline
1,1,1,1,3; 1,1,5 & 21 &0.357\\ \hline
1,1,1,1,3; 3,4 & 15&0.318\\ \hline
1,1,1,1,3; 2,5 & 11 &0.281\\ \hline
1,1,1,1,3; 1,6 & 5& 0.189\\ \hline
1,1,1,1,3; 7 & 1 & 0\\ \hline
1,2,2,2; 1,2,2,2 & 120&0.562 \\ \hline
1,2,2,2; 1,1,2,3 & 87&0.524 \\ \hline
1,2,2,2; 1,1,1,4 & 51& 0.461\\ \hline
1,2,2,2; 2,2,3 & 51 &0.461\\ \hline
1,2,2,2; 1,3,3 & 37&0.424\\ \hline
1,2,2,2; 1,2,4 & 30 &0.399\\ \hline
1,2,2,2; 1,1,5 & 15 &0.318\\ \hline
1,2,2,2; 3,4 & 13 &0.301\\ \hline
1,2,2,2; 2,5 & 9 &0.258\\ \hline
1,2,2,2; 1,6 & 4 &0.163\\ \hline
1,2,2,2; 7 & 1& 0\\ \hline
1,1,2,3; 1,1,2,3 & 67& 0.493\\ \hline
1,1,2,3; 1,1,1,4 &43&0.441\\ \hline
1,1,2,3; 2,2,3 & 39&0.430 \\ \hline
1,1,2,3; 1,3,3 &30 &0.399\\ \hline
1,1,2,3; 1,2,4 & 25 &0.378\\ \hline
1,1,2,3; 1,1,5 &14&0.310\\ \hline
1,1,2,3; 3,4 & 11 &0.281\\ \hline
1,1,2,3; 2,5 & 8&0.244 \\ \hline
1,1,2,3; 1,6 & 4& 0.163 \\ \hline
1,1,2,3; 7 & 1 & 0 \\ \hline
1,1,1,4; 1,1,1,4 & 34 &0.414\\ \hline
1,1,1,4; 2,2,3 & 25 &0.378\\ \hline
1,1,1,4; 1,3,3 &20 &0.351\\ \hline
1,1,1,4; 1,2,4 & 19 &0.345\\ \hline
1,1,1,4; 1,1,5 & 13&0.301 \\ \hline
1,1,1,4; 3,4 & 8 &0.244\\ \hline
1,1,1,4; 2,5 & 7 &0.228\\ \hline
1,1,1,4; 1,6 & 4&0.163\\ \hline
1,1,1,4; 7 & 1& 0 \\ \hline

\end{tabular}

\begin{tabular}{|l|l|l|}
\hline
Messages & Coeff& Deg\\ \hline
2,2,3; 2,2,3 & 25& 0.378\\ \hline
2,2,3; 1,3,3 &19 &0.345\\ \hline
2,2,3; 1,2,4 & 16 &0.325\\ \hline
2,2,3; 1,1,5 &9 &0.258\\ \hline
2,2,3; 3,4 & 3 &0.129\\ \hline
2,2,3; 2,5 &6&0.210\\ \hline
2,2,3; 1,6 & 3 &0.129\\ \hline
2,2,3; 7 & 1 & 0 \\ \hline
1,3,3; 1,3,3 & 15 &0.318\\ \hline
1,3,3; 1,2,4 & 13 &0.301\\ \hline
1,3,3; 1,1,5 & 8 &0.244\\ \hline
1,3,3; 3,4 & 7 &0.228\\ \hline
1,3,3; 2,5 & 5&0.189\\ \hline
1,3,3; 1,6 & 3 &0.129\\ \hline
1,3,3; 7 & 1 & 0\\ \hline
1,2,4; 1,2,4 & 12&0.291\\ \hline
1,2,4; 1,1,5 &8 &0.244\\ \hline
1,2,4; 3,4 & 6 &0.210\\ \hline
1,2,4; 2,5 & 5 &0.189\\ \hline
1,2,4; 1,6 & 3&0.129 \\ \hline
1,2,4; 7 & 1 & 0 \\ \hline
1,1,5; 1,1,5 & 7 &0.228 \\ \hline
1,1,5; 3,4 & 4 &0.163\\ \hline
1,1,5; 2,5 & 4 &0.163\\ \hline
1,1,5; 1,6 & 3 &0.129\\ \hline
1,1,5; 7 &   1& 0 \\ \hline
3,4; 3,4 & 4& 0.163\\ \hline
3,4; 2,5 & 3&0.129\\ \hline
3,4; 1,6 &2 &0.081\\ \hline
3,4; 7 & 1 & 0 \\ \hline
2,5; 2,5 & 3 &0.129\\ \hline
2,5; 1,6 & 2& 0.081 \\ \hline
2,5; 7 & 1&0\\ \hline
1,6; 1,6 & 2 &0.081\\ \hline
1,6; 7 & 1 & 0\\ \hline
7; 7 & 1 & 0\\ \hline
\end{tabular}

\pagebreak

\section{Cardinality of equivalence classes}
\label{appendixB}

We can outline a procedure for calculating the cardinality of the equivalence classes as required for the metric of Gierlichs {\em et al.}  

We know the size of both the sender sets and the receiver sets.  Suppose there are $p$ sender sets and $q$ receiver sets. Each receiver set will consist of  elements from the sender sets in some combination.  

\begin{enumerate}
\item  First determine all integer partitions of each of the receiver sets.  This can be done, for example, with the algorithm from \cite{ivan}.  The integer partitions will count the number of possible ways the sender sets can be partitioned in the receiver set.  For example, if the receiver set has $3$ elements then the partitions are: $1,1,1$; $2,1$; $3$.  We could have one element from each sender set, two elements from one of the sender sets and another element from another sender set, or all elements from a single sender set.  However, the problem here is that determining this properly is probably equivalent to explicitly determining the equivalence classes to begin with and this laborious procedure is something we avoid through the power of our symmetric function method.  We want to find all sequences $t_{1i}, t_{2i},\ldots, t_{ki}$ such that $\sum_a t_{ai}= r_i$.  But we want for the same set of $t$'s to find all sequences $t_{j1}, t_{j2},\ldots, t_{j \ell}$ such that $\sum_b t_{jb}= s_j$.  This could be done either by exhaustive search or through an integer programming approach. 
\item  Since we have determined how many ways each sender set could be partitioned so that $t_{ji}$ elements go to receiver set $i$ for $1\leq i\leq k$, we can construct a series of multinomial coefficients of the form $\frac{s_j !}{t_{j1}!t_{j2}!\ldots t_{j\ell}!}$.   
\item Now the elements within each receiver set of size $r_i$ can be arranged in $r_i!$ ways.
\end{enumerate}

The cardinality of receiver set $i$ is $\prod_j \frac{s_j !}{t_{j1}!t_{j2}!\ldots t_{j\ell}!} r_i!$.   

The cardinality of the equivalence class is $\prod_i \prod_j \frac{s_j !}{t_{j1}!t_{j2}!\ldots t_{j\ell}!} r_i!$ .  


\begin{thebibliography}{ll}
\bibitem{BPS00} O. Berthold, A. Pfitzmann, R. Standtke, The disadvantages of free MIX routes and how to overcome them, in {\em Designing Privacy Enhancing Technologies}, LNCS 2009, 2000, pp. 30--45.
\bibitem{C} D. Chaum, Untraceable electronic mail, return addresses, and digital pseudonyms, {\em Communications of the ACM}24 (1981), 84--88.
\bibitem{D} G. Danezis, Statistical disclosure attacks, in D. Gritzalis, S. De Capitani di Vimercati, P. Samarati, S.K. Katsikas, eds, {\em SEC}, vol 250 of {\em IFIP Conference Proceedings}, Kluwer, 2003, 421--426.

\bibitem{DDT} G. Danezis, C. Diaz, C. Troncoso, Two-sided statistical disclosure attack, {\em Privacy Enhancing Technologies}, 2007.
\bibitem{DSCP02} C. Diaz, S. Seys, J. Claessens, B. Preneel, Towards measuring anonymity," in {\em Proceedings of Privacy Enhancing Technolgies Workshop (PET 2002)}, LNCS 2482, 2002.
\bibitem{ESY07}  M. Edman, F. Sivrikaya, B. Yenner, A Combinatorial Approach to Measuring Anonymity, In {\em Intelligence and Security Informatics}, 2007, pp 356--363.  
\bibitem{Franz} M. Franz, B. Meyer, A. Pashalidis, Attacking unlinkability: the importance of context, Privacy Enhancing Technologies Workshop (PET), 2007.
\bibitem{GTDPV08} B. Gierlichs, C. Troncoso, C. Diaz, B. Preneel, I. Verbauwhede, Revisiting a combinatorial approach towards measuring anonymity, {\em WPES 08}, 2008.
\bibitem{GouldenJackson} I.P. Goulden and D.M. Jackson, {\em Combinatorial Enumeration}, New York: Wiley, 1983.  Dover Reprints 2004.
\bibitem{KAP} D. Kesdogan, D. Agrawal, S. Penz, Limits of anonymity in open environments, in F.A.P. Petitcolas, ed, {\em Information Hiding} LNCS 2578, Springer-Verlag 2002, 53-69.
\bibitem{MacMahon} P.A. MacMahon, {\em Combinatory Analysis}, Vol I, Cambridge University Press, 1915, Reprinted AMS Chelsea, Providence 2001.
\bibitem{RR98} M Reiter and A. Rubin, Crowds: Anonymity for web transactions, {\em ACM Transactions on Information and System Security}, 1 (1998), 66--92.
\bibitem{SD02} A. Serjantov and G. Danezis, Towards an information theoretic measure for anonymity, {\em Proceedings of Privacy Enhancing Technologies Workshop (PET 2002)}, LNCS 2482, 2002.
\bibitem{Stanley} R.P. Stanley, {\em Enumerative Combinatorics}, Vol. II, Cambridge: Cambridge University Press, 1999. 
\bibitem{Stembridge} J.R. Stembridge, The Home Page for SF, posets, and coxeter/weyl, http://www.math.lsa.umich.edu/\~jrs/maple.html, accessed October 6, 2008.
\bibitem{TGPV08} C. Troncoso, B. Gierlichs, B. Preneel, I. Verbauwhede, Perfect Matching Disclosure Attacks, {\em Proceedings of Privacy Enhancing Technologies Workshop (PET 2008)}.
\bibitem{Wilf} H.S. Wilf, {\em Generatingfunctionology}, 3rd ed., A.K. Peters, 2006.
\bibitem{ivan}
A. Zoghbiu and I. Stojmenovic,
\emph{Fast Algorithms For Generating Integer Partitions},
Intern. J. Computer Math., Vol.70. pp. 319---332.

\end{thebibliography}
\end{document}